\documentclass[11pt]{article}
\usepackage{a4wide} 
\usepackage{ae}
\usepackage{graphicx}
\usepackage{amsmath}
\usepackage{amssymb}
\usepackage{amsthm}
\usepackage{amsfonts}
\usepackage{latexsym}
\usepackage{verbatim}
\usepackage{verbatim}
\usepackage{url}
\usepackage[normalem]{ulem}
\advance \topmargin by -\headheight
\advance \topmargin by -\headsep
\evensidemargin \oddsidemargin
\newtheorem{thrm}{Theorem}[section]

\newtheorem{lmm}[thrm]{Lemma}
\newtheorem{crllr}[thrm]{Corollary}
\newtheorem{prpstn}[thrm]{Proposition}

\theoremstyle{definition}
\newtheorem{dfntn}[thrm]{Definition}
\newtheorem{rmrk}[thrm]{Remark}

\newtheorem*{notation}{Notation}

\theoremstyle{remark}
\newtheorem*{note}{Note}

\newcommand{\cA}{{\mathcal A}}

\newcommand{\cP}{{\mathcal P}}

\newcommand{\cC}{{\mathcal C}}

\newcommand{\NN}{{\mathbb N}}

\newcommand{\empt}{\varepsilon}

\newcommand{\rev}{\widetilde}

\newcommand{\bbf}{\mathbf{f}}

\newcommand{\bw}{\mathbf{w}}
\newcommand{\bx}{\mathbf{x}}

\newcommand{\bs}{\mathbf{s}}

\def\outdeg{\mbox{deg}^{+}}
\def\indeg{\mbox{deg}^{-}}

\numberwithin{equation}{section}

\author{Michelangelo Bucci\footnotemark[1] \and Alessandro De Luca\footnotemark[2] \and Amy Glen\footnotemark[3] \and Luca Q.~Zamboni\footnotemark[4]}

\title{A connection between palindromic and factor complexity using return words}

\date{Submitted: February 9, 2008; Accepted: March 25, 2008; Revised: April 10, 2008}

\begin{document}

\normalsize

\maketitle 

\hrule
\begin{abstract}  In this paper we prove that for any infinite word $\bw$ whose set of factors is closed under reversal, the following conditions are  equivalent:
\begin{itemize}
\item[(I)] all complete returns to palindromes are palindromes;
\item[(II)] $\cP(n) + \cP(n+1) = \cC(n+1) - \cC(n) + 2$ for all $n$,
\end{itemize}
where $\cP$ (resp.~$\cC$) denotes the {\em palindromic complexity} (resp.~{\em factor complexity})  function of $\bw$, which counts the number of distinct palindromic factors (resp.~factors) of each length in $\bw$.

\medskip

\vspace{0.1cm} 
\end{abstract}
\hrule

\footnotetext[1]{Dipartimento di Matematica e Applicazoni ``R.~Caccioppoli'', Universit\`a degli Studi di Napoli Federico II, Via Cintia, Monte S.~Angelo, I-80126, Napoli, ITALY ({\tt micbucci@unina.it}).}

\footnotetext[2]{Dipartimento di Matematica e Applicazoni ``R.~Caccioppoli'', Universit\`a degli Studi di Napoli Federico II, Via Cintia, Monte S.~Angelo, I-80126, Napoli, ITALY ({\tt alessandro.deluca@unina.it}).}

\footnotetext[3]{{\bf Corresponding author:} LaCIM, Universit\'e du Qu\'ebec \`a Montr\'eal, C.P. 8888, succursale Centre-ville, Montr\'eal, Qu\'ebec, H3C 3P8, CANADA (\texttt{amy.glen@gmail.com}). Supported by CRM, ISM, and LaCIM.}

\footnotetext[4]{Institut Camille Jordan, 
Universit\'e Claude Bernard Lyon 1, 
43 boulevard du 11 novembre 1918, 
69622 Villeurbanne Cedex
FRANCE ({\tt luca@unt.edu}).}

\section{Introduction}

Given an infinite word $\bw$, let $\cP(n)$ (resp.~$\cC(n)$) denote the {\em palindromic complexity} (resp.~{\em factor complexity}) of $\bw$, i.e., the number of distinct palindromic factors (resp.~factors) of $\bw$ of length $n$.  In \cite{jAmBjCdD03pali}, J.-P. Allouche, M.~Baake, J.~Cassaigne, and D.~Damanik established the
following inequality relating the palindromic and factor complexities of a non-ultimately periodic infinite word: 
\[
\cP(n)\leq \frac{16}{n}\thinspace \cC\hspace{-2pt}\left(n+\left\lfloor{\frac n4}\right\rfloor\right) \quad \mbox{for all $n
\in \NN$}.
\]
More recently, using {\em Rauzy graphs}, P.~Bal\'a\v{z}i, Z.~Mas\'akov\'a, and E.~Pelantov\'a~\cite{pBzMeP07fact} proved that for any uniformly recurrent infinite word whose set of factors is 
 closed under reversal,  
\begin{equation} \label{eq:pali-bound}
  \cP(n) + \cP(n+1) \leq \cC(n+1) - \cC(n) + 2 \quad \mbox{for all $n
\in \NN$}.
\end{equation}
They also provided several examples of infinite words for which $\cP(n) + \cP(n+1)$ always reaches the upper bound given in relation~\eqref{eq:pali-bound}. Such infinite words include {\em Arnoux-Rauzy sequences}, {\em complementation-symmetric sequences}, certain words associated with $\beta$-expansions where $\beta$ is a {\em simple Parry number}, and a class of words coding $r$-interval exchange transformations.

In this paper we give a characterization of all infinite words with factors closed under reversal for which the equality $\cP(n) + \cP(n+1) = \cC(n+1) - \cC(n) + 2$ holds for all $n$: these are exactly the infinite words with the property that all `complete returns' to palindromes are palindromes. Given a finite or infinite word $w$ and a factor $u$ of $w$, we say that a factor $r$ of $w$ is a {\em complete return} to $u$ in $w$ if $r$ contains exactly two occurrences of $u$, one as a prefix and one as a suffix.  Return words play an important role in the study of minimal subshifts; see \cite{Du3, Du2, Du1, FMN, cHlZ99desc, Si}.  

Our main theorem is the following:

\begin{thrm} \label{T:main}
 For any infinite word $\bw$ whose set of factors is closed under reversal, the following  conditions are equivalent: 
\begin{itemize}
\item[{\em (I)}] all complete returns to any palindromic factor of $\bw$ are palindromes;
\item[{\em (II)}] $\cP(n) + \cP(n+1) = \cC(n+1) - \cC(n) +2$ for all $n\in \NN$.
\end{itemize}
\end{thrm}

Recently, in \cite{aGjJ07pali}, it was shown that property (I) is equivalent to every factor $u$ of $\bw$ having exactly $|u| + 1$ distinct palindromic factors (including the empty word). Such words are `rich' in palindromes in the sense that they contain the maximum number of different palindromic factors. Indeed, X.~Droubay, J.~Justin, and G.~Pirillo \cite{xDjJgP01epis} observed that any finite word $w$ of {\em length} $|w|$ contains at most $|w| + 1$ distinct palindromes. 

The family of finite and infinite words having property (I) are called {\em rich words} in \cite{aGjJ07pali}. In independent work, P.~Ambro{\v{z}}, C.~Frougny, Z.~Mas{\'a}kov{\'a},  and E.~Pelantov{\'a}~\cite{pAzMePcF06pali}  have considered the same class of words which they call {\em full words}, following earlier work of S.~Brlek, S.~Hamel, M.~Nivat, and C.~Reutenauer in \cite{sBsHmNcR04onth}. 

Rich words encompass the well-known family of {\em episturmian words} originally introduced by X.~Droubay, J.~Justin, and G.~Pirillo in~\cite{xDjJgP01epis} (see Section~\ref{S:remarks} for more details). Another special class of rich words consists of S.~Fischler's sequences with ``abundant palindromic prefixes'', which were introduced and studied in~\cite{sF06pali} in relation to Diophantine approximation (see also~\cite{sF06pali2}). Other examples of rich words that are neither episturmian nor of ``Fischler type'' include:  non-recurrent rich words, like $abbbb\cdots$ and $abaabaaabaaaab\cdots$; the periodic rich infinite words: $(aab^kaabab)(aab^kaabab)\cdots$, with $k \geq 0$; the non-ultimately periodic recurrent rich infinite word $\psi(\bbf)$ where $\bbf = abaababaaba\cdots$ is the {\em Fibonacci word} and $\psi$ is the morphism: $a \mapsto aab^kaabab$, $b \mapsto bab$; and the recurrent, but not uniformly recurrent, rich infinite word generated by the morphism: $a \mapsto aba$, $b\mapsto bb$. (See~\cite{aGjJ07pali} for these examples and more.)

From the work in \cite{xDjJgP01epis, aGjJ07pali}, we have the following equivalences.

\begin{prpstn} \label{P:rich-aGjJ07pali}  
A finite or infinite word $w$ is {\em rich} if equivalently:
\begin{itemize}
\item all complete returns to any palindromic factor of $\bw$ are palindromes;
\item every factor $u$ of $w$ contains $|u|+1$ distinct palindomes;
\item the longest palindromic suffix of any prefix $p$ of $w$ occurs exactly once in $p$.
\end{itemize}
\end{prpstn}

From the perspective of richness, our main theorem can be viewed as a characterization of {\em recurrent rich infinite words} since any rich infinite word is recurrent if and only if its set of factors is closed under reversal (see~\cite{aGjJ07pali} or Remark~\ref{R:rich-recurrence}). Interestingly, the proof of Theorem~\ref{T:main} relies upon another new characterization of rich words (Proposition~\ref{P:v2reverse}), which is useful for establishing the key step, namely that the so-called {\em super reduced Rauzy graph} is a tree. This answers a claim made in the last few lines of \cite{pBzMeP07fact} where it was remarked that the Rauzy graphs of words satisfying equality~(II) must have a very special form. 

After some preliminary definitions and results in the next section, Section 3 is devoted to the proof of Theorem~\ref{T:main} and some interesting consequences are proved in Section 4.

\section{Preliminaries}

\subsection{Notation and terminology}

In this paper, all words are taken over a finite {\em alphabet} $\cA$, i.e., a finite non-empty set of symbols called {\em letters}. A {\em finite word} over $\cA$ is a finite sequence of letters from $\cA$. The {\em empty word} $\empt$ is the empty sequence. A (right) \emph{infinite word} $\bx$ is a sequence indexed by $\NN_+$ with values in $\cA$, i.e., $\bx = x_1x_2x_3\cdots$ with each $x_i \in \cA$. For easier reading, infinite words are hereafter typed in boldface to distinguish them from finite words. 

Given a finite word $w = x_1x_2\cdots x_m$ (where each $x_i$ is a letter), the  \emph{length} of $w$, denoted by $|w|$, is equal to $m$. By convention, the empy word is the unique word of length $0$. We denote by $\tilde w$ the {\em reversal} of $w$, given by $\tilde w = x_m \cdots x_2x_1$. If $w = \tilde w$, then $w$ is called a {\em palindrome}.

 A finite word $z$ is a {\em factor} of a finite or infinite word $w$ if $w = uzv$ for some words $u$, $v$. In the special case $u = \empt$ (resp.~$v = \empt$), we call $z$ a {\em prefix} (resp.~{\em suffix}) of $w$. If $u \ne \empt$ and $v \ne \empt$, then we say that $z$ is an {\em interior factor} of $w =uzv$. Moreover, $z$ is said to be a {\em central factor} of $w$ if $|u| = |v|$. We say that $z$ is {\em unioccurrent} in $w$ if $z$ occurs exactly once in $w$.  For any finite or infinite word $w$, the set of all factors of  $w$ is denoted by $F(w)$ and we denote by $F_n(w)$ the set of all factors of $w$ of length $n$, i.e., $F_n(w) := F(w) \cap \cA^n$ (where $|w|\geq n$ if $w$ is finite). We say that $F(w)$ is {\em closed under reversal} if for any $u \in F(w)$, $\tilde u \in F(w)$.

 A factor of an infinite word $\bw$ is \emph{recurrent} in $\bw$ if it occurs infinitely often in $\bw$, and $\bw$ itself is said to be \emph{recurrent} if all of its factors are recurrent in it. Furthermore, $\bw$ is \emph{uniformly recurrent} if any factor of $\bw$ occurs infinitely many times in $\bw$ with bounded gaps. 

\begin{rmrk} \label{R:rich-recurrence}
A noteworthy fact (proved in \cite{aGjJ07pali}) is that a rich infinite word is recurrent if and only if its set of factors is closed under reversal.  
 \end{rmrk}
  
  More generally, we have the following well-known result:
  
\begin{prpstn}[folklore] \label{P:folklore}
If $\bw$ is an infinite word with $F(\bw)$ closed under reversal, then $\bw$ is recurrent.
\end{prpstn}
\begin{proof}
Consider some occurrence of a factor $u$ in $\bw$ and let $v$ be a prefix of $\bw$ containing $u$. As $F(\bw)$ is closed under reversal, $\tilde v \in F(\bw)$. Thus, if $v$ is long enough, there is an occurrence of  $\tilde u$ strictly on the right of this particular occurrence of $u$ in $\bw$. Similarly $u$ occurs on the right of this $\tilde u$ and thus $u$ is recurrent in $\bw$. 
\end{proof}

 \subsection{Key results}

We now prove two useful results, the first being a new characterization of rich words.

\begin{prpstn} \label{P:v2reverse} A finite or infinite word $w$ is rich if and only if, for each factor $v \in F(w)$, any factor of $w$ beginning with $v$ and ending with $\tilde v$ and not containing $v$ or $\tilde v$ as an interior factor is a palindrome.
 \end{prpstn}
 \begin{proof}  ONLY IF: Consider any factor $v \in F(w)$ and let $u$ be a factor of $w$ beginning with $v$ and ending with $\tilde v$ and not containing $v$ or $\tilde v$ as an interior factor. If $v$ is a palindrome, then either $u = v = \rev v$ (in which case $u$ is clearly a palindrome), or $u$ is a complete return to $v$ in $w$, and hence $u$ is (again) a palindrome by Proposition~\ref{P:rich-aGjJ07pali}. Now assume that $v$ is not a palindrome.

Suppose by way of contradiction that $u$ is not a palindrome and let $p$ be the longest palindromic suffix of $u$ (which is unioccurrent in $u$ by richness). Then $|p| < |u|$ as $u$ is not a palindrome. If $|p| > |v|$, then $\tilde v$ is a proper suffix of $p$, and hence $v$ is a proper prefix of $p$. But then $v$ is an interior factor of $u$, a contradiction. On the other hand, if $|p| \leq |v|$, then $|p| \ne |v|$ and $p$ is a proper suffix of $\tilde v$ (as $\tilde v$ is not a palindrome), and hence $p$ is a proper prefix of $v$. Thus $p$ is both a prefix and a suffix of $u$; in particular $p$ is not unioccurrent in $u$, a contradiction. \medskip

IF: The given conditions tell us that any complete return to a palindromic factor $v$ ($= \tilde v$) of $w$  is a palindrome. Hence $w$ is rich by Proposition~\ref{P:rich-aGjJ07pali}.
\end{proof}

\begin{prpstn} \label{P:v2v}
 Suppose $w$ is a rich word. Then, for any non-palindromic factor $v$ of $w$, $\tilde v$ is a unioccurrent factor of any complete return to $v$ in $w$.
\end{prpstn}
\begin{proof} Let $r$ be a complete return to $v$ in $w$ and let $p$ be the longest palindromic suffix of $r$. Then $|p| > |v|$; otherwise, if $|p| \leq |v|$, then $p$ would occur at least twice in $r$ (as a suffix of each of the two occurrences of $v$ in $r$), which is impossible as $r$ is rich. Thus $v$ is a proper suffix of $p$, and hence $\tilde v$ is a proper prefix of $p$. So $\tilde v$ is clearly an interior factor of $r$.

It remains to show that $\tilde v$ is unioccurrent in $r$. Arguing by contradiction,  we suppose that $\tilde {v}$ occurs more than once in $r$. Then a complete return $r'$ to $\tilde v$ occurs as a proper factor of $r$. Using the same reasoning as above, $v$ is an interior factor of $r'$, and hence an interior factor of $r$, contradicting the fact that $r$ is a complete return to $v$. Thus $\tilde v$ is unioccurent in $r$.
\end{proof}

\begin{note}
The above proposition tells us that for any factor $v$ of a rich word $w$, occurrences of $v$ and $\tilde v$ alternate in $w$.
\end{note}

\section{Proof of Theorem~\ref{T:main}}

Following the method of Bal\'a\v{z}i {\em et al.}~\cite{pBzMeP07fact}, a key tool for the proof of our main theorem is the notion of  a {\em Rauzy graph}, defined as follows. Given an infinite word $\bw$, the {\em Rauzy graph of order $n$} for $\bw$, denoted by $\Gamma_n(\bw)$, is the directed graph with set of vertices $F_n(\bw)$ and set of edges $F_{n+1}(\bw)$ such that an edge $e \in F_{n+1}(\bw)$  starts at  vertex $v$ and ends at a vertex $v'$ if and only if $v$ is a prefix of $e$ and $v'$ is a suffix of $e$. For a vertex $v$, the {\em out-degree} of $v$ (denoted by $\outdeg(v)$) is the number of distinct edges leaving $v$,    and the {\em in-degree} of $v$ (denoted by $\indeg(v)$) is the number of distinct edges entering $v$. More precisely: 
\[
\outdeg(v) = \sharp\{x \in \cA \mid vx \in F_{n+1}(\bw)\} \quad \mbox{and} \quad \indeg(v) = \sharp\{x \in \cA \mid xv \in F_{n+1}(\bw)\}.
\] 
We observe that, for all $n \in \NN$, 
\[
  \sum_{v \in F_{n}(\bw)} \outdeg(v) = \sharp F_{n+1}(\bw) = \sum_{v \in F_n(\bw)} \indeg(v).
\]  
(Note that $\sharp F_{n+1}(\bw) = \cC(n+1)$.) Hence
\begin{equation} \label{eq:C-diff}
 \cC(n+1) - \cC(n)  = \sum_{v\in F_{n}(\bw)} (\outdeg(v) -1) = \sum_{v\in F_{n}(\bw)} (\indeg(v) -1) .
\end{equation}
It is therefore easy to see that a factor $v \in F_{n}(\bw)$ positively contributes to $\cC(n+1) - \cC(n)$ if and only if $\outdeg(v) \geq 2$, i.e., if and only if there exist at least two distinct letters $a$, $b$ such that $va$, $vb \in F_{n+1}(\bw)$, in which case $v$ is said to be a {\em right-special} factor of $\bw$. Similarly, a factor $v \in F_{n}(\bw)$ is said to be a {\em left-special} factor of $\bw$ if there exist at least two distinct letters $a$, $b$ such that  $av$, $bv \in F_{n+1}(\bw)$. A factor of $\bw$ is said to be {\em special} if it is either left-special or right-special (not necessarily both). With this terminology, if we let $S_n(\bw)$ denote the set of special factors of $\bw$ of length $n$, then formula~\eqref{eq:C-diff} may be expressed as:
\begin{equation}  \label{eq:C-diff2}
 \cC(n+1) - \cC(n)  = \sum_{v \in S_n(\bw)} (\outdeg(v) -1) \quad \mbox{for all $n \in \NN$}.
\end{equation}

Using similar terminology to that in~\cite{pBzMeP07fact}, a directed path $P$ in the Rauzy graph $\Gamma_n(\bw)$ is said to be a {\em simple path of order $n$} if it begins with a special factor $v$ and ends with a special factor $v'$  and contains no other special factors, i.e., $P$ is a directed path of the form $vv'$ or $vz_1\cdots z_kv'$ where each $z_i$ is a non-special factor of length~$n$.  
A special factor $v \in S_n(\bw)$ is called a {\em trivial} simple path of order $n$. 

In what follows, we use the following terminology for paths. Hereafter, ``path'' should be taken to mean ``directed path''.

\begin{dfntn} \label{D:paths} Suppose $\bw$ is an infinite word and let $P= v\cdots v'$ be a path in $\Gamma_n(\bw)$.
\begin{itemize}
\item The first vertex $v$ (resp.~last vertex $v'$)  is called the {\em initial vertex} (resp.~{\em terminal vertex}) of~$P$. 
\item A vertex of $P$ that is neither an initial vertex nor a terminal vertex of $P$ is called an {\em interior vertex} of $P$. 
\item $P$ is said to be a {\em non-trivial path} if it consists of at least two distinct vertices. 
\item The {\em reversal} $\tilde P$ of the path $P$ is the path obtained from $P$ be reversing all edge labels (and arrows) and all labels of vertices. 
\item We say that $P$ is {\em palindromic} (or that $P$ is {\em invariant under reversal}) if $P = \tilde P$. 
\end{itemize}
\end{dfntn}
\begin{note} Given a path $P$ in $\Gamma_n(\bw)$, the reversal of $P$ does not necessarily exist in $\Gamma_n(\bw)$. 
\end{note}

Suppose $P = w_1w_2\cdots w_k$ is a non-trivial path in $\Gamma_n(\bw)$, and for each $i$ with $1\leq i \leq k$, let $a_i$ and $b_i$ denote the respective first and last letters of $w_i$. Then, by the definition of $\Gamma_n(\bw)$, we have $w_1b_2\cdots b_k = a_1\cdots a_{k-1}w_k$. We call this word the {\em label} of the path $P$, denoted by $\ell_P$. Note that the $i$-th shift of $\ell_P := w_1b_2\cdots b_k$ begins with $w_{i+1}$ for all $i$ with $1\leq i \leq k-1$.

For our purposes, it is convenient to consider the {\em reduced Rauzy graph of order $n$}, denoted by $\Gamma_n'(\bw)$, which is the directed graph obtained from $\Gamma_n(\bw)$ by replacing each {\em simple path} $P= w_1w_2\cdots w_{k-1}w_k$ with a directed edge $w_1 \rightarrow w_k$ labelled by $\ell_P$. Thus the set of vertices of $\Gamma_n'(\bw)$ is $S_n(\bw)$. For example, consider the (rich) {\em Fibonacci word}:
\[
  \bbf = abaababaabaababaababaabaababaabaababaababaabaababaaba \cdots
\]
which is generated by the {\em Fibonacci morphism} $\varphi: a\mapsto ab, b \mapsto a$. The reduced Rauzy graph $\Gamma_2'(\bbf)$ consists of the two (special) vertices: $ab$, $ba$ and three directed edges: $ab \to ba$, $ba \to ba$, $ba \to ab$ with respective labels: $aba$, $baab$, $bab$.


\begin{lmm} \label{L:rich-label} Let $\bw$ be a rich infinite word and suppose $P = w_1w_2\cdots w_k$ is a non-trivial path  in $\Gamma_n(\bw)$ with $k\geq 2$. Then the label $\ell_P = w_1b_2\cdots b_{k}$ is a rich word.
\end{lmm}
\begin{proof} We proceed by induction on the number of vertices $k$ in $P$. The lemma is clearly true for $k=2$ since $\ell_P = w_1b_2$ is a factor of $\bw$ of length $n+1$. Now suppose $k\geq 3$ and assume that the label of any path consisting of $k-1$ vertices is rich. Consider any path consisting of $k$ vertices, namely $P=w_1w_2\cdots w_k$, and suppose by way of contradiction that its label $\ell_P = w_1b_2\cdots b_k$ is not rich. Then the longest palindromic prefix $p$ of $\ell_P$ occurs more than once in $\ell_P$. Hence there exists a complete return $r$ to $p$ which is a prefix of $\ell_P$. It follows that $r = \ell_P$, otherwise $r$ would be a factor of the prefix $u := w_1b_2\cdots b_{k-1}$ of $\ell_P$, and hence a palindrome since $u$ is rich by the induction hypothesis. But this  contradicts the maximality of the palindromic prefix $p$. So $\ell_P$ is a non-palindromic complete return to $p$. Let $q$ be the longest palindromic prefix of $u$ (which is unioccurrent in $u$ by richness). If $|p| > |q|$, then $q$ is a proper prefix of $p$, and hence $q$ occurs more than twice in $u$, a contradiction. On the other hand, if $|p| \leq |q|$, then $p$ is a prefix of $q$, and hence $p$ is an interior factor of $\ell_P$ (occurring as a suffix of $q$), a contradiction. Thus $\ell_P$ is rich, as required.
\end{proof}

The proof of Theorem~\ref{T:main} relies upon the following extensions of Propositions~\ref{P:v2reverse}--\ref{P:v2v} to paths.  
 
\begin{lmm} {\em (Analogue of Proposition~\ref{P:v2reverse}.)} \label{L:v2reverse-paths} Suppose $\bw$ is a rich infinite word and let $v$ be any factor of $\bw$ of length $n$. If $P=v\cdots\tilde v$ is a path from $v$ to $\tilde v$ in $\Gamma_n(\bw)$ that does not contain $v$ or $\tilde v$ as an interior vertex, then $P$ is palindromic. This property also holds for paths in $\Gamma_n'(\bw)$.
\end{lmm}
\begin{proof} We first observe that if $P$ consists of a single vertex, then $P = v = \tilde v$, and hence $P$ is palindromic. Now suppose $P$ is a non-trivial path. If $P=v\tilde v$, then $P$ is clearly palindromic. So suppose $P = vz_1\cdots z_k\tilde v$ where the $z_i$ are factors of $\bw$ of length $n$. By definition, the label $\ell_P = vb_1\cdots b_kb_{k+1}$ begins with $v$ and ends with $\tilde v$ and contains neither $v$ nor $\tilde v$ as an interior factor (otherwise $P$ would contain $v$ or $\tilde v$ as an interior vertex, which is not possible). Thus, as $\ell_P$ is rich (by Lemma~\ref{L:rich-label}), it follows   that $\ell_P$ is a palindrome by Proposition~\ref{P:v2reverse}; whence $P$ must be invariant under reversal too. It is easy to see that this property is also true for paths in the reduced Rauzy graph $\Gamma_n'(\bw)$.
\end{proof}

\begin{lmm} {\em (Analogue of Proposition~\ref{P:v2v}.)} \label{L:v2v-paths} Suppose $\bw$ is a  rich infinite word and let $v$ be any non-palindromic factor of $\bw$ of length $n$. If $P=v\cdots v$ is a non-trivial path in $\Gamma_n(\bw)$ that does not contain $v$ as an interior vertex, then $P$ passes through $\tilde v$ exactly once. This property also holds for paths in $\Gamma_n'(\bw)$.
\end{lmm}
\begin{note}
Of particular usefulness is the fact that any path from $v$ to $v$ must pass through $\tilde v$.
\end{note}
\begin{proof} Let us write $P = vz_1\cdots z_kv$ where the $z_i$ are factors of $\bw$ of length $n$. By definition, the label $\ell_P = vb_1\cdots b_kb_{k+1}$ contains exactly two occurrences of $v$, one as a prefix and one as a suffix  (otherwise, if $\ell_P$ contained $v$ as an interior factor, then $v$ would be an interior vertex of $P$, which is not possible). Thus, as $\ell_P$ is rich (by Lemma~\ref{L:rich-label}), it follows that  $\tilde v$ is a unioccurrent (interior) factor of $\ell_P$ by Proposition~\ref{P:v2v}; whence $P$ passes through $\tilde v$ exactly once. It is easy to see that this property is also true for paths in the reduced Rauzy graph $\Gamma_n'(\bw)$.
\end{proof}


\subsection{(I) implies (II)}

Suppose $\bw$ is an infinite word with $F(\bw)$ closed under reversal and satisfying property (I). Then $\bw$ is recurrent by Proposition~\ref{P:folklore} (i.e., $\bw$ is a recurrent rich infinite word). Moreover, recurrence implies that for all $n$, the Rauzy graph $\Gamma_n(\bw)$ is strongly connected, i.e., there exists a directed path from any vertex $v$ to every other vertex $v'$ in $\Gamma_n(\bw)$.

Fix  $n \in \NN$ and let us now consider the {\em super reduced Rauzy graph of order $n$}, denoted by $\Gamma_n''(\bw)$, whose set of vertices consists of all $[v] := \{v, \tilde v\}$ where $v$ is any special factor of length $n$. Any two distinct vertices $[v]$, $[w]$ (with $v \not\in \{w,\tilde w\}$) are joined by an undirected edge with label $[\ell_P] := \{\ell_P, \ell_{\tilde P}\}$ if $P$ or $\tilde P$ is a simple path beginning with $v$ or $\tilde v$ and ending with $w$ or $\tilde w$.  For example, in the case of the Fibonacci word, $\Gamma_2''(\bbf)$ consists of only one vertex: $[ab]$. In general, the super reduced Rauzy graph consists of more than one vertex and may contain multiple edges between vertices.

Suppose $\Gamma_n''(\bw)$ consists of $s$ vertices; namely $[v_i]$, $i = 1, \ldots$, $s$. Since $\Gamma_n(\bw)$ is strongly connected (by recurrence), $\Gamma_n''(\bw)$ is connected; thus it contains at least $s - 1$ edges. 

Now, from Lemma~\ref{L:v2reverse-paths}, we know that if $v$ is a special factor, any simple path from $v$ to $\tilde v$ is palindromic (i.e., invariant under reversal). Moreover, by closure under reversal, if there exists a simple path $P$ from a special factor ${v}$ to a special factor $w$, with $v \not\in \{w,\tilde w\}$, then there is also a simple path from $\tilde w$ to $\tilde v$ (namely, the reversal of the path $P$). Neither of these simple paths is palindromic. 

We thus deduce that there exist {\em at least} $2(s-1)$ non-trivial simple paths in the Rauzy graph $\Gamma_n(\bw)$ that are non-palindromic (i.e., not invariant under reversal). In fact, we will show that there are exactly $2(s-1)$ non-trivial simple paths of order $n$ that are non-palindromic. Indeed, if this true then, as each  palindromic factor of length $n$ or $n+1$ is a central factor of a (unique) palindromic simple path of order $n$, we have:
\begin{equation} \label{eq:1}
  \cP(n) + \cP(n+1) = \sum_{v \in S_n(\bw)} \outdeg(v) - 2(s-1) + p
\end{equation}
where, on the right hand side, the first summand is the total number of non-trivial simple paths, the second summand is the number of non-trivial simple paths that are non-palindromic, and $p$ is the number of special palindromes of length $n$ (i.e., the number of trivial simple paths of order $n$ that are palindromic). By observing that the number of special factors of length $n$ is $2s - p$, we can simplify equation~$\eqref{eq:1}$ to obtain the required equality~(II) as follows:
\begin{eqnarray*}
\cP(n) + \cP(n+1) &=& \sum_{v \in S_n(\bw)} \outdeg(v) - (2s - p) + 2 \\
                         &=& \sum_{v \in S_n(\bw)} (\outdeg(v) - 1) + 2 \\
                         &=& \cC(n+1) - \cC(n) + 2 \qquad \mbox{(by \eqref{eq:C-diff2})}.
\end{eqnarray*}

We observe, in particular, that any infinite word $\bw$ with $F(\bw)$ closed under reversal satisfies equality~(II) if and only if any simple path between a special factor and its reversal is palindromic, and for each $n$, there are exactly $2(s-1)$ non-trivial simple paths of order $n$ that are non-palindromic. The latter condition says that, for all $n$, the super reduced Rauzy graph $\Gamma_n''(\bw)$ contains exactly $s-1$ edges (with each edge corresponding to a simple path and its reversal), and hence $\Gamma_n''(\bw)$ is a tree as it contains $s$ vertices, $s-1$ edges, and must be connected by the recurrence of $\bw$ (which follows from Proposition~\ref{P:folklore}). More formally:

\begin{prpstn} \label{P:equality-satisfied} An infinite word $\bw$ with $F(\bw)$ closed under reversal  satisfies equality~{\em (II)} if and only if the following conditions hold:
\begin{itemize}
\item[$1)$] any simple path between a special factor and its reversal is palindromic; 
\item[$2)$] the super reduced Rauzy graph $\Gamma_n''(\bw)$  is a tree for all $n$.  
\end{itemize}
\end{prpstn}
\begin{proof}  Suppose $\bw$ is an infinite word with $F(\bw)$ closed under reversal. Then $\bw$ is recurrent by Proposition~\ref{P:folklore}. We have already shown that conditions $1)$ and $2)$ imply that $\bw$ satisfies equality~(II). Conversely, if at least one of conditions $1)$ and $2)$ does not hold, then $\cP(n) + \cP(n+1) < \cC(n+1) - \cC(n) + 2$ (by the arguments preceding this proposition), i.e., $\bw$ does not  satisfies equality~(II).
\end{proof}

To complete the proof of ``(I) $\Rightarrow$ (II)'', it remains to show that any recurrent rich infinite word $\bw$ satisfies condition 2) of Proposition~\ref{P:equality-satisfied}, since we have already shown that condition 1) holds for any such $\bw$ (using Lemma~\ref{L:v2reverse-paths}). The proof of the fact that $\bw$ satisfies condition~2) uses the following two lemmas (Lemmas~\ref{L:paths}--\ref{L:path-exists}). 

\begin{notation} Given two distinct special factors $v$, $w$ of the same length $n$, we write $v \not\rightarrow w$ if there does not exist a directed edge from $v$ to $w$ in the reduced Rauzy graph $\Gamma_n'(\bw)$ (i.e., if there does not exist a simple path from $v$ to $w$).
\end{notation}

\begin{lmm} \label{L:paths}
Suppose $\bw$ is a recurrent rich infinite word and let $v$, $w$ be two distinct special factors of $\bw$ of the same length with $v \not\in \{w,\tilde w\}$. If there exists a simple path $P$ from $v$ to $w$, then $P$ is unique and there also exists a unique simple path from $\tilde w$ to $\tilde v$ (namely, the reversal of $P$). Moreover:  
\begin{itemize}
\item[i)] $v \not\rightarrow \tilde w$, and hence $w \not\rightarrow \tilde v $ (unless $w$ is a palindrome);
\item[ii)] $\tilde w \not\rightarrow v$, and hence $\tilde v \not\rightarrow w$ (unless $v$ is a palindrome);
\item[iii)] $w \not\rightarrow v$, and hence $\tilde v \not\rightarrow \tilde w$ (unless $v$ and $w$ are both  palindromes).
\end{itemize}
\end{lmm}
\begin{proof} By closure under reversal (Remark~\ref{R:rich-recurrence}), if there exists a simple path $P$ from $v$ to $w$, then the reversal of $P$ is a simple path from  $\tilde w$ to $\tilde v$ in the Rauzy graph of order $|v|=|w|=n$. To prove the uniqueness of $P$, let us suppose there exist two different simple paths $P_1$, $P_2$ from $v$ to $w$ in the Rauzy graph $\Gamma_n(\bw)$. Then
\begin{eqnarray*}
P_1 = vu_1\cdots u_kw \quad \mbox{and} \quad P_2 = vz_1\cdots z_\ell w \quad \mbox{for some $k$, $l \in \NN$},
\end{eqnarray*}
where $u_1$, \ldots, $u_k$, $z_1$, \ldots, $z_\ell$ are non-special factors of $\bw$ of length $n$ and $u_i \ne z_i$ for some $i$. Note that either $P_1$ or $P_2$ (not both) may be of the form $vw$.

To keep the rest of the proof as simple as possible, we assume hereafter that neither $v$ nor $w$ is a palindrome; the arguments are similar, and in fact easier, in the cases when either $v$ or $w$ (or both) is a palindrome.

Consider a path $Q$ of minimal length beginning with $P_1$ and ending with~$P_2$ (in the Rauzy graph $\Gamma_n(\bw)$):
\[
Q = P_1\cdots P_2 = vu_1\cdots u_k\underbrace{w \cdots v}_{ Q_1}z_1\cdots z_\ell w.
\]
First we observe that $Q$ contains $\tilde v$ since any path from $v$ to itself must pass through $\tilde v$, by Lemma~\ref{L:v2v-paths}. Moreover, the left-most $\tilde v$ in $Q$ must occur in the subpath $ Q_1$ (since $\tilde v$ is not equal to any of the non-special factors $u_i$, $z_j$ and $\tilde v \ne w$). Therefore
\[
  Q = \underbrace{vu_1\cdots u_kw\cdots \tilde v}_{Q_2} \cdots vz_1\cdots z_\ell w
 \]
 where the subpath $ Q_2$ ends with the left-most $\tilde v$ in the path $Q$. By Lemma~\ref{L:v2v-paths}, $Q_2$ is a path from $v$ to $\tilde v$ that does not contain $v$ or $\tilde v$ as an interior vertex. Thus, by Lemma~\ref{L:v2reverse-paths}, $Q_2$ is palindromic, and hence $ Q_2$ ends with the reversal of the path $P_1$ since it begins with $P_1$. More explicitly:
 \[
 Q = \underbrace{vu_1\cdots u_kw}_{P_1}\cdots \overbrace{\underbrace{\tilde w\tilde u_k \cdots \tilde u_1\tilde v}_{\rev P_1} \cdots \underbrace{vz_1\cdots z_\ell w}_{P_2}}^{Q_3}.
 \]
 We distinguish two cases. \medskip
 
\noindent{\it Case $1$}: If the subpath $Q_3$ contains $w$ as a terminal vertex only, then $\tilde w$ is not an interior vertex of  $Q_3$ by Lemma~\ref{L:v2v-paths}, and hence $Q_3$ is palindromic by Lemma~\ref{L:v2reverse-paths}. It follows that $k = \ell$ and $z_i = u_{i}$ for all $i = 1, \ldots, k$. Thus $P_1 = P_2$; a contradiction. \medskip

\noindent{\it Case $2$}: If the subpath $Q_3$ contains $w$ as an interior vertex, then $Q_3$ first passes through $w$ after taking the path $\rev P_1$ (at the beginning) and before taking the path $P_2$ (at the end). Hence, by Lemma~\ref{L:v2reverse-paths}, $Q_3$ begins with a palindromic path from $\tilde w$ to $w$ that begins with $\rev P_1$ and hence ends with $P_1$. But then $Q$ passes through the path $P_1$ at least twice before taking the path $P_2$, contradicting the fact that $Q$ is a path of minimal length beginning with $P_1$ and ending with $P_2$. \medskip

Both cases lead to a contradiction; thus the simple path $P$ from $v$ to $w$ is unique (and its reversal $\rev P$ is the unique simple path from $\tilde w$ to $\tilde v$). It remains to show that conditions $i)$--$iii)$ hold. As $ii)$ is symmetric to $i)$, we prove only that $i)$ and $iii)$ are satisfied. By what precedes, it suffices to consider paths in the reduced Rauzy graph $\Gamma_n'(\bw)$.

\medskip

\noindent $i)$: Arguing by contradiction, let us suppose that there exists a (unique) simple path from $v$ to $\tilde w$, i.e., there exists a directed edge from $v$ to $\tilde w$ in the reduced Rauzy graph $\Gamma_n'(\bw)$. Then (from above) we know that there also exists a directed edge from $w$ to $\tilde v$. Consider a shortest path $Q$ in the reduced Rauzy graph $\Gamma_n'(\bw)$ beginning with $v\tilde w$ and ending with $vw$. 
By Lemma~\ref{L:v2v-paths}, any path from $v$ to itself passes through $\tilde v$, so we may write
\[
Q = \underbrace{v\tilde w \cdots \tilde v}_{Q_1} \cdots vw,
\]
where the subpath $Q_1$ ends with the left-most $\tilde v$ in the path $Q$. By Lemmas~\ref{L:v2reverse-paths}--\ref{L:v2v-paths}, the path $Q_1 = v\tilde w \cdots \tilde v$ is palindromic, and hence it ends with $w\tilde v$. So we have $Q = v\tilde w \cdots w \tilde v \cdots vw$; moreover, by Lemma~\ref{L:v2v-paths}, $\tilde w$ must occur between the last two $w$'s shown here. In particular, 
\[
 Q = v\tilde w \cdots \underbrace{w\tilde v \cdots \tilde w}_{Q_2} \cdots vw
\]
where the subpath $Q_2$ contains $\tilde w$ as a terminal vertex only. Thus, by Lemmas~\ref{L:v2reverse-paths}--\ref{L:v2v-paths}, the path $Q_2 = w\tilde v \cdots \tilde w$ is palindromic, and hence it ends with $v \tilde w$. But then $Q$ ends with a shorter path of the form $v \tilde w \cdots v w$, contradicting the fact that $Q$ is a path of minimal length beginning with $v\tilde w$ and ending with $vw$. \medskip

\noindent $iii)$: Again, the proof proceeds by contradiction. Suppose there exists a (unique) simple path from $w$ to $v$. Consider a shortest path $Z$ in the reduced Rauzy graph $\Gamma_n'(\bw)$ beginning with $wv$ and ending with $vw$. 
By Lemma~\ref{L:v2v-paths}, the path $Z$ must pass through $\tilde w$; thus
\[
Z = \underbrace{wv \cdots \tilde w}_{Z_1} \cdots vw.
\]
where the subpath $Z_1$ ends with the left-most $\tilde w$ in the path $Z$. Now it follows from Lemmas~\ref{L:v2reverse-paths}--\ref{L:v2v-paths} that the subpath $Z_1$ is palindromic, and hence $Z_1$ must end with $\tilde v\tilde w$. So we may write
\[
Z = wv \cdots \underbrace{\tilde v \tilde w \cdots v}_{Z_2}w.
\]
If the subpath $Z_2$ contains $v$ as a terminal vertex only, then neither $v$ nor $\tilde v$ is an interior vertex of $Z_2$ by Lemma~\ref{L:v2v-paths}. Thus $Z_2$ is palindromic  by Lemma~\ref{L:v2reverse-paths}, and hence $Z_2$ ends with $wv$. But then the path $Z$ ends with the path $w v w$, which is impossible by Lemma~\ref{L:v2v-paths}. Thus, the subpath $Z_2$ must pass through $v$ at an earlier point, and hence we have $Z_2 = \tilde v \tilde w\cdots v\cdots v$. In particular, the path $Z_2$ begins with a palindromic subpath of the form $\tilde v \tilde w \cdots wv$, by Lemma~\ref{L:v2reverse-paths}.  
But then the path $Z$ ends with a shorter path from $wv$ to $vw$, contradicting the minimality of  $Z$.
\end{proof}

\begin{notation} For a finite word $v$, let $v^{\epsilon}$ represent either $v$ or $\tilde v$ and set $v^{-\epsilon} := \rev{v^\epsilon}$.
\end{notation}
\begin{lmm} \label{L:path-exists} Let $\bw$ be a recurrent rich infinite word. For fixed $n \in \NN_+$, suppose the super reduced Rauzy graph $\Gamma_n''(\bw)$ contains at least three distinct  vertices: $[v_1]$, $[v_2]$, \ldots, $[v_s]$, $s \geq 3$. Then, for each $k$ with $3 \leq k \leq s$, the reduced Rauzy graph $\Gamma_n'(\bw)$ contains a path from $v_1$ to $v_k^{\epsilon_k}$ of the form:
\[
v_1v_2^{\epsilon_2}\cdots v_2v_3^{\epsilon_3} \cdots v_{k-2}v_{k-1} ^{\epsilon_{k-1}}\cdots v_{k-1} v_{k}^{\epsilon_k},
\]
where for all $i = 2, \ldots, k-1$, the subpath $v_{i}^{\epsilon_{i}} \cdots v_{i}$ (which may consist of only the single vertex $v_i^{\epsilon_i}$) does not contain $v_j$, $\tilde v_j$ for all $j$ with $1\leq j \leq k$, $j \ne i$.  
\end{lmm}
\begin{proof} We use induction on $k$ and employ similar reasoning to the proof of Lemma~\ref{L:paths}.

First consider the case $k = 3$. Recurrence implies that $\Gamma_n'(\bw)$ is connected, so we may assume without loss of generality that $\Gamma_n'(\bw)$ contains a directed edge from $v_1$ to $v_2^{\epsilon_2}$, a directed edge from $v_2$ to $v_3^{\epsilon_3}$, and a path from $v_2^{\epsilon_2}$ to $v_2$. That is,  $\Gamma_n'(\bw)$ contains a path beginning with $v_1v_2^{\epsilon_2}$ and ending with $v_2v_3^{\epsilon_3}$.  Consider such a path of minimal length:
\[
Q = v_1v_2^{\epsilon_2} \cdots v_2v_3^{\epsilon_3}.
\]
To prove the claim for $k = 3$, we show that none of the special factors $v_1$, $\tilde v_1$, $v_3$, $\tilde v_3$ are interior vertices of $Q$. If $v_2^{\epsilon_2} = v_2$, then $Q=v_1v_2v_3^{\epsilon_3}$ (by  minimality) and we are done. So let us assume that  $v_2^{\epsilon_2} = \tilde v_2 \ne v_2$.

 Observe that if $v_1$ is an interior vertex of $Q$, then $\tilde v_1$ must be an interior vertex of $Q$ since any path from $v_1$ to itself must contain $\tilde v_1$, by Lemma~\ref{L:v2v-paths}. Similarly, if $v_3^{\epsilon_3}$ is an interior vertex of $Q$, then $v_3^{-\epsilon_3}$ is an interior vertex of $Q$. Therefore it suffices to show that $\tilde v_1$ and  $v_3^{-\epsilon_3}$ are not interior vertices of $Q$. We prove this fact only for $\tilde v_1$ as the proof is similar for $v_3^{-\epsilon_3}$. 

Arguing by contradiction, suppose $\tilde v_1$ is an interior vertex of  $Q$. Then $Q$ begins with a palindromic path from $v_1\tilde v_2$ to $\tilde v_1$ (by Lemmas~\ref{L:v2reverse-paths}--\ref{L:v2v-paths}), and this palindromic path clearly ends with $v_2\tilde v_1$. Hence 
\[ 
Q = v_1\tilde v_2 \cdots \underbrace{v_2\tilde v_1\cdots v_2}_{Q'}v_3^{\epsilon_3}
\]
where the subpath $Q'$ begins with a palindromic path from $v_2\tilde v_1$ to $\tilde v_2$ (by Lemmas~\ref{L:v2reverse-paths}--\ref{L:v2v-paths}), and this palindromic path clearly ends with $v_1\tilde v_2$. But then the path $Q$ ends with a shorter path from $v_1\tilde v_2$ to $v_2v_3^{\epsilon_3}$, contradicting the minimality of $Q$. Thus the lemma holds for $k = 3$.

Now suppose $4 \leq k\leq s$ and assume the claim holds for $k-1$. Since $\Gamma_n'(\bw)$ is connected, it contains a path beginning with $v_1v_2^{\epsilon_2}\cdots v_2v_3^{\epsilon_3} \cdots v_{k-2} v_{k-1}^{\epsilon_{k-1}}$ and ending with $v_{k-1}v_k^{\epsilon_{k}}$ (where the former path satisfies the conditions of the lemma). Consider such a path of minimal length:
\begin{equation}\label{eq:Z}
  Z = \underbrace{v_1v_2^{\epsilon_2}\cdots v_2v_3^{\epsilon_3} \cdots v_{k-2}}_{Z_1} \underbrace{v_{k-1}^{\epsilon_{k-1}} \cdots v_{k-1}}_{Z_2}v_k^{\epsilon_{k}}
\end{equation}
where for all $i = 2, \ldots, k-2$, the subpath $v_{i}^{\epsilon_{i}} \cdots v_{i}$ (which may consist of only the single vertex $v_i^{\epsilon_i}$) does not contain  $v_j$, $\tilde v_j$ for all $j$ with $1\leq j \leq k-1$, $j \ne i$.  To prove the induction step, we show that the path $Z$ satisfies the following two conditions:
\begin{itemize} 
\item[$i)$] the subpath $Z_1$ contains neither $v_k$ nor $\tilde v_k$;
\item[$ii)$] the subpath $Z_2 = v_{k-1}^{\epsilon_{k-1}} \cdots v_{k-1}$ does not contain $v_j$, $\tilde v_j$ for all $j$ with $1\leq j \leq k$, $j \ne k-1$. 
\end{itemize}

First suppose that condition $i)$ is not satisfied, i.e.,  $Z_1$ contains $v_k$ or $\tilde v_k$. Without loss of generality we assume that $v_k$ is the right-most of the vertices $v_k$, $\tilde v_k$ appearing in $Z_1$. 

\medskip
\noindent{\em Case $1$:} Suppose ${v_k^{\epsilon_k}} = v_k \ne \tilde v_k$. Then $Z$ ends with a path from $v_k$ to itself, which must pass through $\tilde v_k$ by Lemma~\ref{L:v2v-paths}; moreover, $\tilde v_k$ must be an interior vertex of $Z_2$ (by the choice of $v_k$). Thus, by Lemmas~\ref{L:v2reverse-paths}--\ref{L:v2v-paths}, $Z_2v_k$ (and hence $Z$) ends with a palindromic path from $\tilde v_k$ to $v_{k-1}v_{k}$. 
Hence $Z_2$ contains $\tilde v_k\tilde v_{k-1}$, and we have:
\[
  Z_2v_k= \underbrace{v_{k-1}^{\epsilon_{k-1}} \cdots \tilde v_k\tilde v_{k-1}}_{Z_3} \cdots v_{k-1}v_k
  \]
where the subpath $Z_3$ ends with a palindromic path from $v_{k-1}$ to $\tilde v_k\tilde v_{k-1}$ (by Lemmas~\ref{L:v2reverse-paths}--\ref{L:v2v-paths}); thus $Z_3$ contains $v_{k-1}v_{k}$. But then $Z$ begins with a shorter path from $Z_1$ to $v_{k-1}v_{k}^{\epsilon_k}$, contradicting the minimality of $Z$. 

\medskip
\noindent{\em Case $2$:} Suppose $v_k^{\epsilon_k} = \tilde v_k$. Then the path $Z$ ($= Z_1Z_2\tilde v_k$)  ends with a path of the form:
\[
 Z_4 = v_k \underbrace{\cdots\cdots}_{\mathrm{no}~v_k, \tilde v_k}Z_2\tilde v_k.
 \]
 If $v_k$ or $\tilde v_k$ is an interior vertex of $Z_2$, then  we reach a contradiction using the same arguments as in Case 1. 
 On the other hand, if neither $v_k$ nor $\tilde v_k$ is an interior vertex of $Z_2$, then  $Z_4$ is palindromic by Lemma~\ref{L:v2reverse-paths}. So the path $Z_4$ begins with $v_k\tilde v_{k-1}$ since it ends with $v_{k-1}\tilde v_k$. But then $\tilde v_{k-1}$ is an interior vertex of $Z_1$, a contradiction.  \medskip

Thus the path $Z$ satisfies condition $i)$. In proving this fact, we have also shown that $v_k$, $\tilde v_k$ are not interior vertices of $Z_2$. It remains to show that the subpath $Z_2$ does not contain $v_j$, $\tilde v_j$  for all $j$ with $1\leq j \leq k-2$ (and hence $Z$ satisfies condition $ii)$). We prove only that $Z_2$ does not contain $\tilde v_1$ or $\tilde v_1$ since the proof is similar when considering other $v_j$, $\tilde v_j$. 

Suppose on the contrary that $Z_2$ contains $v_1$ or $\tilde v_1$. Then, by Lemmas~\ref{L:v2reverse-paths}--\ref{L:v2v-paths}, $Z$ begins with a palindromic path  from $v_1$ to $\tilde v_1$, and this palindromic path  begins with $Y = Z_1v_{k-1}^{\epsilon_{k-1}}$ (and hence ends with $\tilde Y$) by the conditions on $Z$ under the induction hypothesis. More explicitly, we have:
 \[
 Z = \overbrace{\underbrace{v_1v_2^{\epsilon_2} \cdots v_{k-2} v_{k-1}^{\epsilon_{k-1}}}_{Y}\cdots\underbrace{v_{k-1}^{-\epsilon_{k-1}}\tilde v_{k-2} \cdots v_2^{-\epsilon_2}\tilde v_1}_{\tilde Y}}^{\mathrm{palindromic}}\underbrace{\cdots v_{k-1}v_{k}^{\epsilon_k}}_{Z_5}.
\] 
 Hence, as $v_{k-1}$ and $\tilde v_{k-1}$ are not interior vertices of $Y$ (by the induction hypothesis), the subpath $\tilde Y Z_5$ begins with a palindromic path from $v_{k-1}^{-\epsilon_k}$ to $v_{k-1}^{\epsilon_{k-1}}$, and this palindromic path begins with $\tilde Y$ (and hence ends with~$Y$), by Lemmas~\ref{L:v2reverse-paths}--\ref{L:v2v-paths}. But then $Z$ ends with a shorter path from $Y$ to $v_{k-1}v_k^{\epsilon_k}$, contradicting the minimality of~$Z$. 
 
We conclude that the subpath $Z_2 = v_{k-1}^{\epsilon_{k-1}} \cdots v_{k-1}$ does not contain $v_j$, $\tilde v_j$ for all $j$ with $1\leq j \leq k$, $j \ne i$ (i.e., the path $Z$ satisfies condition $ii)$),   and the proof is thus complete.
\end{proof}

\begin{lmm} \label{L:tree} Suppose $\bw$ is a recurrent rich infinite word. Then the super reduced Rauzy graph $\Gamma_n''(\bw)$ is a tree  for all $n \in \NN_+$.
\end{lmm}
\begin{proof} First recall that for all $n$, $\Gamma_n''(\bw)$ is connected (by the recurrence property of $\bw$). Moreover, Lemma~\ref{L:paths} tells us that if two distinct vertices in $\Gamma_n''(\bw)$ are joined by an edge, then this edge is unique (and corresponds to a simple path and its reversal). It remains to show that $\Gamma_n''(\bw)$ does not contain any {\em cycle} (i.e., does not contain a chain linking a vertex with itself).

Suppose on the contrary  that  $\Gamma_n''(\bw)$ contains a cycle for some $n$. Then $\Gamma_n''(\bw)$ must contain at least three distinct vertices:  
$[v_1]$, $[v_2]$, \ldots, $[v_s]$, $s \geq 3$, 
and a cycle of the following form:
\begin{equation} \label{eq:circuit}
\mbox{$[v_1]$---$[v_2]$---~$\cdots$~---$[v_k]$---$[v_1]$ ~~ for some $k$ with $3 \leq k \leq s$}.
\end{equation}
%
We thus deduce from Lemma~\ref{L:path-exists} that the reduced Rauzy graph $\Gamma_n'(\bw)$ contains a path from $v_1$ to $v_1^{\epsilon_1}$ of the form:
\[
P =  v_1v_2^{\epsilon_2}\cdots v_2v_3^{\epsilon_3} \cdots v_{k-2} v_{k-1}^{\epsilon_{k-1}} \cdots v_{k-1}v_k^{\epsilon_{k}}\cdots v_kv_1^{\epsilon_1},
\]
where for all $i = 2, \ldots, k$, the subpath $v_{i}^{\epsilon_{i}} \cdots v_{i}$ (which may consist of only the single vertex $v_i^{\epsilon_i}$) does not contain  $v_j$, $\tilde v_j$ for all $j$ with $1\leq j \leq k$, $j \ne i$.  (Note that $P$ corresponds to the cycle given in \eqref{eq:circuit}.)

First suppose that $v_1$ is a palindrome. In this case, as neither $v_1$ nor $\tilde v_1$ is an interior vertex of $P$, it must be a palindromic path by Lemma~\ref{L:v2reverse-paths}. But then $v_k = v_2^{-\epsilon_2}$, a contradiction (as $k\geq 3$).

Now suppose that $v_1$ is not a palindrome. If $v_1^{\epsilon_1} = \tilde v_1$, then we deduce (as above, using Lemma~\ref{L:v2reverse-paths}) that the path $P$ must be palindromic, yielding a contradiction. On the other hand, if $v_1^{\epsilon} = v_1$, then, by Lemma~\ref{L:v2v-paths}, the path $P$ must pass through $\tilde v_1$, a contradiction.

Thus $\Gamma_n''(\bw)$ is a tree.
\end{proof}

This concludes our proof of the ``(I) $\Rightarrow$ (II)'' part of Theorem~\ref{T:main}.

\subsection{(II) implies (I)}

Conversely, suppose $\bw$ is an infinite word with $F(\bw)$ closed under reversal and satisfying equality~(II). Then $\bw$ satisfies conditions $1)$ and $2)$ of Proposition~\ref{P:equality-satisfied}.  

Now, arguing by contradiction, suppose $\bw$ does not satisfy property (I) (i.e., $\bw$ is not rich). Then there exists a palindromic factor $p$ that has a non-palindromic complete return $u$ in $\bw$; in particular, we have $u = pqavb\tilde qp$ for some words $q$, $v$ (possibly empty) and letters $a$, $b$, with $a \ne b$. So the words $pqa$, $b\tilde q p$ and their reversals $a\tilde q p$, $pqb$ are factors of $\bw$. Thus $pq$ (resp.~$\tilde q p$) is a right-special (resp.~left-special) factor of $\bw$. Hence, if $u$ does not contain any other special factors, then $u$ forms the label of a non-palindromic simple path beginning with $pq$ and ending with $\tilde qp$. But this contradicts condition~1) of Proposition~\ref{P:equality-satisfied}.  Therefore $u$ must contain other special factors of length $n := |pq|$, besides $pq$ and $\tilde qp$. In particular,  $u$ begins with the label of a simple path of order $n$ beginning with $pq$ and ending with another special factor $s_1$ of length $n$. Similarly, $u$ ends with the label of a  simple path of order $n$ beginning with a special factor $s_2$ of length $n$ and ending with $\tilde q p$. Moreover, since $u$ is a complete return to $p$, neither  $s_1$ nor $s_2$ is equal to $pq$ or $\tilde q p$ (otherwise $p$ occurs as an interior factor of $u$). Thus, in the super reduced Rauzy graph $\Gamma_n''(\bw)$, there is an edge between the  vertex $[pq]$ and each of the vertices $[s_1]$ and $[s_2]$. In particular, there exists a path of the form:  $\mbox{$[s_1]$---$[pq]$---$[s_2]$}$.  Furthermore, as $u$ contains a factor that begins with $s_1$ and ends with $s_2$ and contains no occurrence of $pq$ or $\tilde q p$, there also exists a chain (or possibly just an edge) linking $[s_1]$ and $[s_2]$ that does not contain the vertex $[pq]$. Thus, if $\{s_1, \tilde s_1\} \ne \{s_2, \tilde s_2\}$, then we see that $\Gamma_n''(\bw)$ 
contains a cycle, contradicting condition~2) of Proposition~\ref{P:equality-satisfied}. On the other hand, if $\{s_1, \tilde s_1\} = \{s_2, \tilde s_2\}$, then there are at least two edges joining the vertices $[s_1]$ and $[pq]$. Indeed, there exists a simple path $P_1$ from $pq$ to $s_1$ and there also exists a simple path $P_2$ either from $s_1$ to $\tilde q p$ or from $\tilde s_1$ to $\tilde q p$. By closure under reversal, the reversals $\tilde P_1$, $\tilde P_2$ of the respective simple paths $P_1$, $P_2$ also exist. Moreover, none of these four simple paths coincide. Certainly, $P_1 \ne P_2$,  $P_1 \ne \tilde P_1$, and $P_2 \ne \tilde P_2$ as neither $s_1$ nor $\tilde s_1$ is equal to $pq$ or $\tilde q p$, and $P_1 \ne \tilde P_2$ as the second vertex in $P_1$ ends with the letter $a$, whereas the second vertex in the path $\tilde P_2$ ends with the letter $b \ne a$. So $\Gamma_n''(\bw)$ is not a tree, contradicting condition~$2)$ of Proposition~\ref{P:equality-satisfied}. This concludes our proof of Theorem~\ref{T:main}. \qed

\section{A few consequences and remarks} \label{S:remarks}

From Theorem~\ref{T:main}, we easily deduce that property (I) is equivalent to equality (II) for any uniformly recurrent infinite word. Indeed, equality (II) implies the existence of arbitrarily long palindromes since $\cP(n) + \cP(n+1)\geq 2$ for all $n$, so together with uniform recurrence one can readily show that factors are closed under reversal; hence property (I) holds by Theorem~\ref{T:main}. Conversely, richness (property (I)) together with uniform recurrence implies closure under reversal by Remark~\ref{R:rich-recurrence}, and hence equality (II) holds. 

\medskip
\noindent{\bf Question:} {\it In the statement of Theorem~\ref{T:main}, can the hypothesis of factors being  closed under reversal be replaced by the weaker hypothesis of recurrence?} 
\medskip

As above, it follows directly from Theorem~\ref{T:main} and Remark~\ref{R:rich-recurrence} that for any recurrent infinite word $\bw$, if $\bw$ satisfies property (I) (i.e., if $\bw$ is rich, and hence has factors closed under reversal), then equality (II) holds. However, to prove the converse using our methods, one would need to know that any recurrent infinite word satisfying equality (II) has factors closed under reversal. We could not find a proof of this claim nor could we find a counter-example. Let us point out  that whilst uniform recurrence and the existence of arbitrarily long palindromes imply closure under reversal, this is not true in the case of recurrence only. For instance, consider the following infinite word: 
\[
\bs = bca^2bca^3bca^2bca^4bca^2bca^3bca^2bca^5bc\cdots, 
\]
which is the limit as $n$ goes to infinity of the sequence $(s_n)_{n\geq1}$ of finite words defined by:
\[
 s_1 = bc \quad \mbox{and} \quad s_n = s_{n-1}a^ns_{n-1} \quad \mbox{for $n>1$}.
 \]
This infinite word is clearly recurrent (but not uniformly recurrent) and contains arbitrarily long palindromes, but its set of factors is not closed under reversal. (Note that $\bs$ is not rich and does not satisfy equality (II).) If one could show that recurrence together with equality~(II) implies arbitrarily long palindromic prefixes, this would be enough to prove that factors are closed under reversal.

In the context of finite words $w$, the hypothesis of factors being closed under reversal can be replaced by the requirement that $w$ is a palindrome. Indeed, all we really need is the super reduced Rauzy graph to be connected, which is true for palindromes. 

\newpage
\begin{thrm} \label{T:main-finite}
For any palindrome $w$, the following properties are equivalent:
\begin{itemize}
\item[$i)$] $w$ contains $|w| + 1$ distinct palindromes;
\item[$ii)$] all complete returns to palindromes in $w$ are palindromes;
\item[$iii)$] $\cP(i)+\cP(i+1)=\cC(i+1)-\cC(i)+2$ for all $i$ with $0\leq i \leq |w|$. \qed
\end{itemize} 
\end{thrm}

We now prove  two easy consequences of Theorem~\ref{T:main}.

\begin{crllr} Suppose $\bw$ is a recurrent rich infinite word. Then the following properties hold.
\begin{itemize}
\item[$i)$] $\bw$ is (purely) periodic if and only if $\cP(n) + \cP(n+1) = 2$ for some $n$.
\item[$ii)$]  $(\cP(n))_{n\geq 1}$ is eventually periodic with period $2$ if and only if there exist non-negative integers $K$, $L$, $N$ such that $\cC(n) = Kn + L$ for all $n \geq N$. 
\end{itemize}
\end{crllr}
\begin{proof} Suppose $\bw$ is a recurrent rich infinite word. Then $\cP(n) + \cP(n+1) = \cC(n+1) - \cC(n) + 2$ for all $n$, by Theorem~\ref{T:main} and Remark~\ref{R:rich-recurrence}.  \medskip

\noindent $i)$: If $\cP(n) + \cP(n+1) = 2$ for some $n$, then $\cC(n+1) = \cC(n)$, and hence $\bw$ is eventually periodic; in particular, $\bw$ must be (purely) periodic as it is recurrent. Conversely, if $\bw$ is periodic, then  $\cC(n+1) = \cC(n)$ for some $n$, and hence $\cP(n) + \cP(n+1) = 2$. \medskip

\noindent $ii)$:  The condition on $\cC(n)$ implies that for all $n \geq N$, $\cC(n+1) - \cC(n) = K$, and hence  $\cP(n) + \cP(n+1) = K +2 = \cP(n+1) + \cP(n+2)$. Thus $\cP(n) = \cP(n+2)$ for all $n\geq N$. Conversely, suppose $(\cP(n))_{n\geq 1}$ is eventually periodic with period $2$. Then there exists a non-negative integer $N$ such that $\cP(n) = \cP(n+2)$ for all $n \geq N$. Hence, for all $n\geq N$, $\cP(n) + \cP(n+1) = \cC(n+1) - \cC(n) + 2 = \cP(n+1) + \cP(n+2) = M \geq 2$. Therefore $\cC(n+1) - \cC(n) = M - 2$ for all $n \geq N$. 
\end{proof}

\begin{rmrk} Item $ii)$ of the above corollary can be compared with a result of J.~Cassaigne \cite{jC96spec}, who proved that if $\cC(n)$ has linear growth, then $\cC(n+1) - \cC(n)$ is bounded. 
\end{rmrk}

\begin{rmrk} In~\cite{pBzMeP07fact}, Bala\v{z}i {\it et al.}~remarked: ``According to our knowledge, all known examples of infinite words which satisfy the equality $\cP(n) + \cP(n +1) = \cC(n+1) - \cC(n) + 2$ for all $n\in \NN$ have sublinear factor complexity.'' Actually, there do exist recurrent rich infinite words with non-sublinear complexity. For instance, the following example from~\cite{aGjJ07pali}: $abab^2abab^3abab^2abab^4abab^2abab^3abab^2abab^5\cdots$ (which is the fixed point of the morphism: $a\mapsto abab$, $b \mapsto b$) is a recurrent rich infinite word and its complexity $\cC(n)$ grows quadratically with $n$. Another example that was indicated to us by J.~Cassaigne is the fixed point of $ a \mapsto aab$, $b\mapsto b$:  
\[ 
aabaabbaabaabbbaabaabbaabaabbbbaabaabbaabaabbbaabaabbaabaabbbbb\cdots. 
\]
It is a recurrent rich infinite word and its complexity is equivalent to $n^2/2$. More precisely, $\cP(n) + \cP(n+1) - 2 = \cC(n+1) - \cC(n) = n+1 - \sharp\{k>0 ~|~ 2^k+k-2<n\}$.
\end{rmrk}

In \cite{xDjJgP01epis}, X.~Droubay {\em et al.}~showed that the family of {\em episturmian words} (e.g., see \cite{xDjJgP01epis, jJgP02epis, aGjJ07epis}), which includes the well-known {\em Sturmian words}, comprises a special class of uniformly recurrent rich infinite words. Specifically, they proved that if an infinite word $\bw$ is episturmian, then any factor $u$ of $\bw$ contains exactly $|u| + 1$ distinct palindromic factors (see \cite[Cor.~2]{xDjJgP01epis}). An alternative proof of the richness of episturmian words can be found in the paper \cite{vAlZiZ03pali} where the fourth author, together with V.~Anne and I.~Zorca, proved that for episturmian words,  all complete returns to palindromes are palindromes. (A shorter proof of this fact is also given in~\cite{mBaDaDlZ08onso}.) More recently, P.~Bal\'a\v{z}i {\em et al.}~\cite{pBzMeP07fact} showed that all {\em strict} episturmian words (i.e., {\em Arnoux-Rauzy sequences} 
\cite{pAgR91repr, gR82suit}) satisfy $\cP(n) + \cP(n+1) = \cC(n+1) - \cC(n) +2$ for all~$n$. This fact, together with Theorem~\ref{T:main}, provides yet another proof that {\em all} episturmian words are rich (since any factor of an episturmian word is a factor of some strict episturmian word).

Sturmian words are exactly the aperiodic episturmian words over a $2$-letter alphabet. They have complexity $n+1$ for each $n$ and are characterized by their palindromic complexity: any Sturmian word has $\cP(n) = 1$ whenever $n$ is even and $\cP(n) = 2$ whenever $n$ is odd (see~\cite{xDgP99pali}). From these observations, one can readily check that Sturmian words satisfy equality~(II) (and hence they are rich).

We can now say even more: the set of factors of all Sturmian words satisfies equality~(II). To show this, we first recall that F.~Mignosi~\cite{fM91onth} proved that, for any $n\geq 0$, the number $c(n)$ of {\em finite Sturmian words} of length $n$ is given by
\[
 c(n) = 1 + \sum_{i=1}^n(n+1-i)\phi(i),
\]
where $\phi$ is {\em Euler's totient function}. More recently, in \cite{aDaD06comb}, the second author together with A.~de~Luca proved that for any $n\geq 0$, the number $p(n)$ of Sturmian palindromes of length $n$ is given by
\[
 p(n) = 1 + \sum_{i=0}^{\lceil n/2 \rceil -1}\phi(n-2i).
\]
Equivalently, for any $n\geq 0$, 
\[
 p(2n) = 1 + \sum_{i=1}^n\phi(2i) \quad \mbox{and} \quad p(2n+1) = 1 + \sum_{i=0}^n \phi(2i+1).
\]
Thus, for all $n\geq 0$,
\[
  p(2n) + p(2n+1) = 2 + \sum_{i=1}^n\left\{\phi(2i) + \phi(2i+1)\right\} + 2 = \sum_{i=1}^{2n+1}\phi(i) + 2,
\]  
and
\begin{eqnarray*}
  c(2n+1) - c(2n) + 2 &=& \sum_{i=1}^{2n+1}(2n+2-i)\phi(i) - \sum_{i=1}^{2n}(2n+1-i)\phi(i) + 2 \\
                        &=& \phi(2n+1) + \sum_{i=1}^{2n}\phi(i) + 2 \\
                        &=& \sum_{i=1}^{2n+1}\phi(i) + 2
                        = p(2n) + p(2n+1).
\end{eqnarray*}

From this point of view, it would be interesting to count for instance the number of all binary rich words of length $n$ for each $n$.

\bigskip
\noindent {\large\bf Acknowledgements.} The authors would like to thank Jacques Justin for helpful comments and suggestions on a preliminary version of this paper. The first three authors would also like to acknowledge the hospitality of the Department of Mathematics at the University of North Texas where this work was done.


\begin{thebibliography}{99}

\bibitem{jAmBjCdD03pali} J.-P.~Allouche, M.~Baake, J.~Cassaigne, D.~Damanik, Palindrome complexity, {\it Theoret. Comput. Sci.} {292} (2003) 9--31.

\bibitem{pAzMePcF06pali} P.~Ambro{\v{z}}, C.~Frougny, Z.~Mas{\'a}kov{\'a}, E.~Pelantov{\'a},  Palindromic complexity of infinite words associated with 
  simple {P}arry numbers, {\it Ann. Inst. Fourier (Grenoble)} 56 (2006) 2131--2160.
  
  \bibitem{vAlZiZ03pali} V.~Anne, L.Q.~Zamboni, I.~Zorca, Palindromes and
  pseudo-palindromes in episturmian and pseudo-palindromic infinite words, in: {\it Proceedings of the Fifth International Conference on Words} (Montr\'eal, Canada), September 13--17, 2005. {\it Publications du LaCIM} 36 (2005) 91--100. 
  
  \bibitem{pAgR91repr} P.~Arnoux, G.~Rauzy, Repr\'{e}sentation 
g\'{e}om\'{e}trique de suites de complexit\'{e} $2n+1$, 
\emph{Bull. Soc. Math. France} {119} (1991) 199--215.

\bibitem{pBzMeP07fact} P.~Bal\'a\v{z}i, Z.~Mas\'akov\'a, E.~Pelantov\'a, Factor versus
  palindromic complexity of uniformly recurrent infinite words, {\em Theoret. 
  Comput. Sci.} {380} (2007) 266--275. 
  

  \bibitem{sBsHmNcR04onth} S.~Brlek, S.~Hamel, M.~Nivat, C.~Reutenauer,  On the palindromic complexity of infinite words, {\it Internat. J. Found. Comput. Sci.} { 15} (2004) 293--306.
  
  \bibitem{mBaDaDlZ08onso} M.~Bucci, A.~de Luca, A.~De~Luca, L.Q.~Zamboni, On some problems related to palindromic closure, {\it Theoret. Inform. Appl.} (in press), doi:10.1051/ita:2007064.
  
\bibitem{jC96spec} J.~Cassaigne, Special factors of sequences with linear subword complexity, in: {\em Developments in Language Theory II}, World Scientific, Singapore, 1996, pp.~25--34.
  
\bibitem{aDaD06comb}
A.~de~Luca, A. De~Luca, Combinatorial properties of
  {S}turmian palindromes, {\em Internat. J. Found. Comput. 
  Sci.} 17 (2006) 557--573.
  
\bibitem{xDjJgP01epis} X.~Droubay, J.~Justin, G.~Pirillo,
Episturmian words and some constructions of de Luca and Rauzy,
\emph{Theoret. Comput. Sci.} {255} (2001) 539--553.    

\bibitem{xDgP99pali}
X.~Droubay, G.~Pirillo, Palindromes and {S}turmian words, {\em Theoret.
  Comput. Sci.} 223 (1999) 73--85.
  
  
    \bibitem{Du3} F.~Durand, A characterization of substitutive sequences using return words, {\em Discrete Math.} {179} (1998) 89--101.
  
  \bibitem{Du2} F.~Durand, A generalization of Cobham's theorem, {\em Theory Comput. Syst.} {31} (1998) 169--185.
  
  \bibitem{Du1} F.~Durand, Linearly recurrent subshifts have a finite number of non-periodic subshift factors, {\em Ergodic Theory Dynam. Systems} {19} (1999) 953--993.
  
  \bibitem{FMN} S.~Ferenczi, C.~Mauduit, A.~Nogueira, Substitutional dynamical systems: algebraic characterization of eigenvalues, {\em Ann. Sci. \'Ecole Norm. Sup.} {29} (1995) 519--533.
  
  
  \bibitem{sF06pali}
S.~Fischler, Palindromic prefixes and episturmian words,
{\it J. Combin.  Theory Ser.~A} 113 (2006) 1281--1304.

\bibitem{sF06pali2}
S.~Fischler, Palindromic prefixes and diophantine approximation, {\it Monatsh. Math.} 151 (2007) 11--37.
  
\bibitem{aGjJ07epis} A.~Glen, J.~Justin, Episturmian words: a survey, Preprint, 2007, 	arXiv:0801.1655.

\bibitem{aGjJ07pali} A.~Glen, J.~Justin, S.~Widmer, L.Q.~Zamboni, Palindromic richness, {\it European J. Combin.}, to appear, arXiv:0801.1656. 


\bibitem{cHlZ99desc} C.~Holton, L.Q.~Zamboni, Descendants of primitive substitutions, {\em Theory Comput. Syst.} { 32} (1999) 133--157. 



\bibitem{jJgP02epis} J.~Justin, G.~Pirillo, Episturmian words and episturmian morphisms, \emph{Theoret. Comput. Sci.} {276} (2002) 281--313. 

\bibitem{fM91onth} F.~Mignosi, On the number of factors of {S}turmian words, 
  {\em Theoret. Comput. Sci.} 82 (1991) 71--84.
  
\bibitem{gR82suit} G.~Rauzy, Suites \`a termes dans un alphabet fini, in: {\it S\'emin. Th\'eorie des Nombres}, Exp. No. 25, pp.~16, Univ. Bordeaux I, Talence, 1982--1983.

\bibitem{Si} A.~Siegel, Pure discrete spectrum dynamical systems and periodic tiling associated with a substitution, {\em Ann. Inst. Fourier (Grenoble)} {54} (2004) 341--381.

\end{thebibliography}
\end{document}